\documentclass[12pt,reqno]{amsart}

\usepackage[arrow,matrix,curve]{xy}

\usepackage[dvips]{graphicx} 

\usepackage{amssymb, latexsym, amsmath, amscd, array, hyperref
}

\newtheorem{theorem}{Theorem}[section]

\theoremstyle{definition}

\newtheorem{remark}[theorem]{Remark}

\numberwithin{equation}{section}
\newcommand\N {{\mathbb N}} 

\newcommand\R {{\mathbb R}}
\newcommand\Q {{\mathbb Q}}
\newcommand\Z {{\mathbb Z}}

\DeclareMathOperator{\adequal}{\,{}_{\ulcorner\!\urcorner}\,}

\author[P.B.]{Piotr B\l{}aszczyk}\address{P. B\l{}aszczyk, Institute
of Mathematics, Pedagogical University of Cracow,
Poland}\email{pb@up.krakow.pl}

\author[V. K.]{Vladimir Kanovei} \address{V. Kanovei, IPPI, Moscow,
and MIIT, Moscow, Russia}\email{kanovei@googlemail.com}

\author[K. K.]{Karin U. Katz}\address{K. Katz, Department of
Mathematics, Bar Ilan University, Ramat Gan 52900
Israel}\email{katzmik@math.biu.ac.il}

\author[M. K.]{Mikhail G. Katz}\address{M. Katz, Department of
Mathematics, Bar Ilan University, Ramat Gan 52900
Israel}\email{katzmik@macs.biu.ac.il}

\author[T. K.]{Taras Kudryk}\address{T. Kudryk, Department of
Mathematics, Lviv National University, Lviv,
Ukraine}\email{kudryk@mail.lviv.ua}

\author[T. M.]{Thomas Mormann}\address{T. Mormann, Department of Logic
and Philosophy of Science, University of the Basque Country UPV/EHU,
20080 Donostia San Sebastian, Spain} \email{ylxmomot@sf.ehu.es}

\author[D. S.]{David Sherry}\address{D. Sherry, Department of
Philosophy, Northern Arizona University, Flagstaff, AZ 86011,
US}\email{David.Sherry@nau.edu}

\begin{document}

\thispagestyle{empty}

%\huge

\title[Leibnizian calculus and first order logic]{Is Leibnizian
calculus embeddable in first order logic?}

\begin{abstract}
To explore the extent of embeddability of Leibnizian infinitesimal
calculus in first-order logic (FOL) and modern frameworks, we propose
to set aside ontological issues and focus on procedural questions.
This would enable an account of Leibnizian procedures in a framework
limited to FOL with a small number of additional ingredients such as
the relation of infinite proximity.  If, as we argue here, first order
logic is indeed suitable for developing modern proxies for the
inferential moves found in Leibnizian infinitesimal calculus, then
modern infinitesimal frameworks are more appropriate to interpreting
Leibnizian infinitesimal calculus than modern Weierstrassian ones.

\end{abstract}

\maketitle

\tableofcontents

\section{Introduction}

Leibniz famously denied that infinite aggregates can be viewed as
\emph{wholes}, on the grounds that they would lead to a violation of
the principle that the whole is greater than the part.  Yet the
infinitary idea is latent in Leibniz in the form of a distinction
between \emph{assignable} and \emph{inassignable} quantities
\cite[p.~153]{Ch}, and explicit in his comments as to the violation of
Definition~V.4 of Euclid's \emph{Elements} \cite[p.~322]{Le95b}.  This
definition is closely related to what is known since \cite{Sto83} as
the \emph{Archimedean property}, and was translated by Barrow in 1660
as follows:
\begin{quote}
\emph{Those numbers are said to have a ratio betwixt them, which being
multiplied may exceed one the other} \cite{Eu60}.
\end{quote}

Furthermore, Leibniz produced a number of results in infinitesimal
calculus which, nowadays, are expressed most naturally by means of
quantifiers that range over infinite aggregates.  This tension leads
us to examine a possible relationship between Leibnizian infinitesimal
calculus and a modern logical system known as first order logic (FOL).
The precise meaning of the term is clarified in Section~\ref{FOL}.  We
first analyze several Leibnizian examples in Section~\ref{one}.  

This text continues a program of re-evaluation of the history of
infinitesimal mathematics initiated in \cite{KK12}, \cite{B11} and
elsewhere.

\section{Examples from Leibniz}
\label{one}

Let us examine some typical examples from Leibniz's infinitesimal
calculus so as to gauge their relationship to FOL.

\subsection{Series presentation of~$\pi/4$}
\label{s11}

In his \emph{De vera proportione} (1682), Leibniz represented
$\frac{\pi}{4}$ in terms of the infinite series
\[
1-\frac{1}{3}+\frac{1}{5}-\frac{1}{7}+\ldots
\]
This is a remarkable result, but we wish to view it as a result
concerning a specific real number, i.e., a \emph{single case}, and in
this sense involving no quantification, once we add a new function
symbol for a black box procedure~$\square$ called ``evaluation of
convergent series'' (as well as a definition of~$\pi$) (we will say a
few words about the various implementations of~$\square$ in modern
frameworks in Section~\ref{two}).

\subsection{Leibniz convergence criterion for alternating series}
\label{s22}

This refers to an arbitrary alternating series defined by an
alternating sequence with terms of decreasing absolute value tending
to zero, such as the series of Subsection~\ref{s11}, or the series
$\sum_n \frac{(-1)^n}{n}$ determined by the alternating sequence
$\frac{(-1)^n}{n}$.  We will refer to such a sequence as a `Leibniz
sequence' for the purposes of this subsection.  This criterion seems
to be quantifying over sequences (and therefore sets), thus
transcending the FOL framework, but in fact this can be handled easily
by introducing a free variable that can be interpreted later according
to the chosen domain of discourse.

Thus, the criterion fits squarely within the parameters
of~$\text{FOL}+\square$ at level~(3) (see Section~\ref{FOL}).  In
more detail, we are not interested here in arbitrary `Leibniz
sequences' with possibly inassignable terms.  Leibniz only dealt with
sequences with ordinary (assignable) terms, as in the two examples
given above.  Each real sequence is handled in the
framework~$\mathbb{R}\subset{}^\ast\mathbb{R}$ by the transfer
principle, which asserts the validity of each true relation when
interpreted over~${}^\ast\mathbb{R}$.

\subsection{Product rule}
\label{s23}

We have~$\frac{d(uv)}{dx}=\frac{du}{dx}v+\frac{dv}{dx}u$ and it looks
like we need quantification over pairs of functions~$(u,v)$.  Here
again we are only interested in natural extensions of real
functions~$u,v$, which are handled at level~(3) as in the previous
section.

In \cite{Le84}, the product rule is expressed in terms of
differentials as~$d(uv)=udv+vdu$.  In \emph{Cum Prodiisset}
\cite[p.~46-47]{Le01c}, Leibniz presents an alternative justification
of the product rule (see \cite[p.~58]{Bos}).  Here he divides by~$dx$
and argues with differential quotients rather than differentials.
%
%\footnote{\label{arthur1}Leibniz freely inverts his infinitesimals,
%making it difficult to interpret his infinitesimals in terms of modern
%nilsquare ones, as Arthur attempts to do in \cite{Ar} (see also
%footnote~\ref{arthur2}).}
%
Adjusting Leibniz's notation, we obtain an equivalent calculation
\[
\begin{aligned}
\frac{d(uv)}{dx} &= \frac{(u+du)(v+dv)-uv}{dx} \\&=
\frac{udv+vdu+du\,dv}{dx} \\&= \frac{udv+vdu}{dx} + \frac{du\,dv}{dx}
\\&= \frac{udv+vdu}{dx}.
\end{aligned}
\]
Under suitable conditions the term~$\frac{du\,dv}{dx}$ is
infinitesimal, and therefore the last step
\begin{equation}
\label{72}
\frac{udv+vdu}{dx} + \frac{du\,dv}{dx} =
u\,\frac{dv}{dx}+v\,\frac{du}{dx},
\end{equation}
relying on a generalized notion of equality, is legitimized as an
instance of Leibniz's \emph{transcendental law of homogeneity}, which
authorizes one to discard the higher-order terms in an expression
containing infinitesimals of different orders.

\subsection{Law of continuity}

Leibniz proposed a heuristic principle known as the \emph{law of
continuity} to the effect that
\begin{quote}
\ldots{} et il se trouve que les r\`egles du fini r\'eussissent dans
l'infini \ldots{} ; et que vice versa les r\`egles de l'infini
r\'eussissent dans le fini, \ldots{} \cite[p.~93-94]{Le02a},
\end{quote}
cited by \cite[p.~67]{Kn02}, \cite[p.~262]{Ro66},
\cite[p.~145]{Lau92}, and other scholars.

On the face of it, one can find numerous counterexamples to such a
principle.  Thus, finite ordinal number addition is commutative,
whereas for infinite ordinal numbers, the addition is no longer
commutative:~$1 + \omega = \omega \not= \omega+1$.  Thus, the infinite
realm of Cantor's ordinals differs significantly from the finite: in
the finite realm, commutativity rules, whereas in the infinite, it
does not not.  Thus the transfer of properties between these two
realms fails.

Similarly, there are many infinitary frameworks where the law of
continuity fails to hold.  For example, consider the Conway--Alling
surreal framework; see e.g., \cite{Al85}.  Here one can't extend even
such an elementary function as~$\sin(x)$ from~$\R$ to the surreals.
Even more strikingly,~$\sqrt{2}$ turns out to be (sur)rational; see
\cite[chapter~4]{Co01}.  The surnaturals don't satisfy the Peano
Arithmetic.  Therefore transfer from the finite to the infinite domain
fails also for the domain of the surreals.

On the other hand, the combined insight of \cite{He48}, \cite{Lo55},
and \cite{Ro61} was that there does exist an infinitary framework
where the law of continuity can be interpreted in a meaningful
fashion.  This is the~$\R\subset{}^\ast\R$ framework.  While it is not
much of a novelty that many infinitary systems don't obey a law of
continuity/transfer, the novelty is that there is one that does, as
shown by Hewitt, {\L}o{\'s}, and Robinson, in the context of
first-order logic.

Throughout the 18th century, Euler and other mathematicians relied on
a broad interpretation of the law of continuity or, as Cauchy will
call it, the \emph{generality of algebra}.  This involved manipulation
of infinite series as if they were finite sums, and in some cases it
also involved ignoring the fact that the series diverges.  The first
serious challenge to this principle emerged from the study of Fourier
series when new types of functions arose through the summation
thereof.  Specifically, Cauchy rejected the principle of the
\emph{generality of algebra}, and held that a series is only
meaningful in its radius of convergence.  Cauchy's approach was
revolutionary at the time and immediately attracted followers like
Abel.  Cauchy in 1821 was perhaps the first to challenge such a broad
interpretation of the law of continuity, with a possible exception of
Bolzano, whose work dates from only a few years earlier and did not
become widely known until nearly half a century later.  For additional
details on Cauchy see \cite{KK11}, \cite{KK12}, \cite{BK},
\cite{Ba14}.  For Euler see \cite{KKKS} and \cite{Ba16b}.

\subsection{Non-examples: EVT and IVT}

It may be useful to illustrate the scope of the relevant results by
including a negative example.  Concerning results such as the extreme
value theorem (EVT) and the intermediate value theorem (IVT), one
notices that the proofs involve procedures that are not easily encoded
in first order logic.  These 19th century results (due to suitable
combinations of Bolzano, Cauchy, and Weierstrass) arguably fall
outside the scope of Leibnizian calculus, as do infinitesimal
foundations for differential geometry as developed in \cite{NK},
\cite{KKN}.

There are axioms in FOL for a real closed field~$F$ (e.g., real
algebraic numbers, real numbers, hyperreal numbers, Conway numbers).
One of these axioms formalizes the fact that IVT holds for odd degree
polynomials~$F[x]$.  In fact, one needs infinitely many axioms like
$(\forall a,b,c) (\exists x) [x^3+ax^2+bx+c=0]$.  Meanwhile, IVT in
its full form is equivalent to the continuity axiom for the real
numbers \cite{Bl15}.

\section{What does ``first-order'' mean exactly?}
\label{FOL}

The adjective `first-order' as we use it entails limitations on
quantification over sets (as opposed to elements).  Now Leibniz really
did not have much to say about properties of sets in general in the
context of his infinitesimal calculus, and even declared on occasion
that infinite totalities don't exist, as mentioned above.  Note that
Leibniz arguably did exploit second-order logic in areas outside
infinitesimal calculus (see \cite{Len87}, \cite{Len04}) but this will
not be our concern here.  Leibniz famously takes for granted second
order logic in formulating his principle governing the \emph{identity
of indiscernibles}.  While second order logic is possibly part of
Leibniz's metaphysics it is not in any obvious way part of his
infinitesimal calculus.

Once we reach topics like Baire category, measure theory, Lebesgue
integration, and modern functional analysis, quantification over sets
becomes important, but these were not Leibnizian concerns in the kind
of analysis he explored.

In fact, the term ``first order logic" has several meanings.  We can
distinguish three levels at which a number system could have
first-order properties compatible with those of the real numbers.
Note that the real numbers satisfy the axioms of an ordered field as
well as a completeness axiom.
\begin{enumerate}
\item
An \emph{ordered field} obeys those among the usual axioms of the real
number system that can be stated in first-order logic (completeness is
excluded).  For example, the following commutativity axiom
holds:~$(\forall x, y)\;[x+y=y+x]$.
\item
A real closed ordered field has all the first-order properties of the
real number system, regardless of whether such properties are usually
taken as axiomatic, for statements which involve the basic
ordered-field relations~$+, \times$, and~$\leq$. This is a stronger
condition than obeying the ordered-field axioms. More specifically,
one includes additional first-order properties, such as extraction of
roots (e.g., existence of a root for every odd-degree polynomial).
For example, every number must have a cube root:~$(\forall x)(\exists
y)\;[y^3=x]$, or every positive number have a square root: $(\forall
x>0)(\exists y)\; [y^2=x]$.
\item
The system could have all the first-order properties of the real
number system for statements involving arbitrary relations (regardless
of whether those relations can be expressed using~$ +, \times$, and
$\leq$).  For example, there would have to be a sine function that is
well defined for infinitesimal and infinite inputs; the same is true
for every real function.  To do series, one needs a symbol for~$\N$,
so as to define transcendental entities such as~$\pi$ or sine.  We
also introduce function symbols for whatever functions we are
interested in working with; say all elementary functions occurring in
Leibniz as well as their combinations via composition,
differentiation, and integration.
\end{enumerate}

It follows from these examples that the \emph{first order}
qualification is connected with the intended domain of discourse, so
that any quantifier related to objects outside the domain of discourse
is qualified as not a first order one.  It could be added however that
all mathematical objects are, generally speaking, (represented by
suitable) sets from the set theoretic standpoint, and hence all
mathematical quantifiers are first-order with respect to the
background set universe (superstructure).

The point with level (3) is that instead of quantifying over sequences
or functions, we relate to each individual sequence or function, and
make sure that it has an analogue in the extended domain.  Such an
analogue of~$f$ is sometimes referred to as the \emph{natural
extension} of~$f$.  Then we can say something about the extension of
every standard object in our system, e.g., function, without ever
being able to assert anything about \emph{all} functions.  Thus, the
product rule for differentiation is proved for the assortment of
functions chosen in item (3) above.

\begin{remark}
An alternative to the multitude of functional symbols would be to add
a countable list of variables $u,v,w,\ldots$ meant to denote
unspecified functions.  The idea is to avoid quantifying over such
variables, and use them as merely free variables.  Then, for example,
the product rule is the following statement: ``if $u,v$ are
differentiable functions then the Leibniz rule holds for $u$ and
$v$'', with $u,v$ being free variables.
\end{remark}

Note that we use FOL in a different sense from that used in
formalizing Zermelo--Fraenkel set theory (ZFC).

%and involves rather what Robinson describes as Henkin semantics,
%involving the distinction between internal and external sets.

When we seek hyperreal proxies, following the pioneering work of
\cite{He48} and \cite{Ro66}, for Leibniz's procedural moves, the
theory of \emph{real closed fields} at level~(2) is insufficient and
we must rely upon level~(3).

Thus, Leibniz's series of Subsection~\ref{s11} is expressible
in~$\text{FOL}+\square$ at level (3) but FOL level (2) does not
suffice since~$\pi$ is not algebraic.  Similarly, examples in
Subsection~\ref{s22} and Subsection~\ref{s23} need symbols for
unspecified functions which are not available at level (2).

\section{Modern frameworks}
\label{two}

Based on the examples of Section~\ref{one}, we would like to consider
the following question:
\begin{quote}
\emph{Which modern mathematical framework is the most appropriate for
interpreting Leibnizian infinitesimal calculus?}
\end{quote}
The frameworks we would like to consider are
\begin{enumerate}
\item[(A)] a Weierstrassian (or ``epsilontic'') framework in the
context of what has been called since \cite{Sto83} an
\emph{Archimedean} continuum, satisfying Euclid V.4 (see
Section~\ref{one}), namely the real numbers exclusively; and
\item[(B)] a modern framework exploiting infinitesimals such as the
hyperreals, which could be termed a \emph{Bernoullian} continuum since
Johann Bernoulli was the first to exploit infinitesimals (rather than
``exhaustion'' methods) systematically in developing the calculus.%
\footnote{The adjective \emph{non-Archimedean} is used in modern
mathematics to refer to certain modern theories of ordered number
systems properly extending the real numbers, namely various successors
of \cite{Ha07}.  In modern mathematics, this adjective tends to evoke
associations unrelated to 17th century mathematics.  Furthermore,
defining infinitesimal mathematics by a negation, i.e., as
\emph{non-Archimedean}, is a surrender to the
Cantor--Dedekind--Weierstrass (CDW) view.  Meanwhile, true
infinitesimal calculus as practiced by Leibniz, Bernoulli, and others
is the base of reference as far as 17th century mathematics is
concerned.  The CDW system could be referred to as non-Bernoullian,
though the latter term has not yet gained currency.}
\end{enumerate}
The series summation blackbox~$\square$ (see Subsection~\ref{s11}) is
handled differently in A and B.  Framework~A exploits a first-order
``epsilontic'' formulation that works in a complete Archimedean field.
Thus, the convergence of a series~$\sum_i u_i$ to~$L$ would be
expressed as follows:
\[
(\forall\epsilon>0)(\exists n\in\N)(\forall m\in\N)\;\left[m\geq
n\rightarrow \left|\sum_{i=1}^m u_i-L\right|<\epsilon\right].
\]
\cite[Chapter~5]{Is90} sought to interpret Leibnizian infinitesimal
calculus by means of such quantified paraphrases, having apparently
overlooked Leibniz's remarks to the effect that his infinitesimals
violate Euclid V.4 \cite[p.~322]{Le95b}.

Meanwhile, framework~B allows for an alternative interpretation in
terms of the shadow (i.e., the standard part, closely related to
Leibniz's generalized notion of equality) and hyperfinite partial sums
as follows: for each infinite hypernatural~$H$ the partial sum
$\sum_{i=1}^H u_i$ is infinitely close to~$L$, i.e.,
\[
(\forall H) \left[H \text{ infinite } \rightarrow \sum_{i=1}^H
u_i\approx L\right].
\]
This is closer to the historical occurrences of the package~$\square$
as found in Gregory, Leibniz, and Euler, as we argue below.

In pursuing modern interpretations of Leibniz's work, a helpful
distinction is that between ontological and procedural issues.  More
specifically, we seek to sidestep traditional questions concerning the
ontology of mathematical entities such as numbers, and concentrate
instead on the procedures, in line with Quine's comment to the effect
that
\begin{quote}
Arithmetic is, in this sense, all there is to number: there is no
saying absolutely what the numbers are; there is only arithmetic.
\cite[p.~198]{Qu}
\end{quote}
Related comments can be found in \cite{Be65}.  If one could separate
the ``ontological questions" from the rest, then framework~A would be
more appropriate than framework~B for interpreting the classical texts
\emph{if and only if} framework~A provides better proxies for the
procedural moves found in Leibnizian infinitesimal calculus than
framework~B does, and vice versa.

The tempting evidence in favor of the appropriateness of a modern
framework~B for interpreting Leibnizian infinitesimal calculus is the
presence of infinitesimals and infinite numbers in both, as well as
the availability of hyperreal proxies for guiding principles in
Leibniz's work such as the law of continuity as expressed in
\cite{Le01c} and \cite{Le02a} as well as the transcendental law of
homogeneity \cite{Le10b}; see \cite{KS2}, \cite{KS1}, \cite{SK},
\cite{Gu}.  To what extent Leibnizian infinitesimals can be
implemented in differential geometry can be gauged from \cite{NK}.

The question we seek to explore is whether the limitation of working
with first order logic as discussed in Section~\ref{FOL} could
potentially undermine a full implementation of a hyperreal scheme for
Leibnizian infinitesimal calculus.

With this in mind, let us consider Skolem's construction of
nonstandard natural numbers \cite{Sk33}, \cite{Sk34}, \cite{Sk55}; see
\cite[section~3.2]{KKM} for additional references.  It turns out that
one needs many, many nonstandard numbers in order to move
from~$\mathbb{N}$ to~${}^\ast\mathbb{N}$, e.g., in Henkin's countable
model one has
\begin{equation}
\label{21}
{}^\ast\mathbb{N}=\mathbb{N} + (\mathbb{Z}\times\mathbb{Q}).
\end{equation}
Here we use~$\Q$ to indicate that the galaxies are dense, so that
between any pair of galaxies there is another galaxy (a galaxy is the
set of numbers at finite distance from each other).  Meanwhile~$\Z$
indicates that each galaxy other than the original~$\N$ itself is
order-isomorphic to~$\Z$ rather than to~$\N$, because for each
infinite~$H$ the number~$H-1$ is in the same galaxy.

Leibniz arguably did not have such a perspective.  In other words, one
needs to build up a considerable conceptual machinery to emulate
Leibniz's probably rather modest arsenal of procedural moves.  That is
to say, we may be able to emulate all of Leibniz moves in a modern
B-framework, such as Leibniz's infinite quantities, his distinction
between assignable and inassignable quantities, and his transcendental
law of homogeneity.  However, the B-framework also enables us to carry
out many additional moves unknown to Leibniz, for instance those
related to the detailed structure of~${}^\ast\N$ as in \eqref{21}.
Thus, the \emph{difference} between the Leibnizian framework and a
modern infinitesimal B-framework is \emph{large}.

On the other hand, a Weierstrassian A-framework may not cover all the
moves Leibniz may make in his framework LEI, but one might argue that
the \emph{difference} between (A) and LEI is small.  Thus, one may not
necessarily have~$\text{LEI}\subseteq\text{(A)}$, but one might argue
that the difference~$\text{(A)}-\text{LEI}$ is \emph{small}.  This may
be taken as evidence that (A) and LEI are more similar to each other
than (B) and LEI are.  This could affect the assessment of
appropriateness.  Finally, could it be that neither the Weierstrassian
nor the modern infinitesimal account is appropriate to cope with
Leibnizian infinitesimal calculus?

\section{Separating entities from procedures}
\label{three}

What would it mean exactly to separate ontological problems from
procedural problems?  A possible approach is to attempt to account for
Leibniz's procedures in a framework limited to first order logic, with
a small number of additional ingredients such as the relation of
\emph{infinite proximity} and the closely related \emph{shadow
principle} for passing from a finite inassignable quantity to an
assignable one (or from a finite nonstandard number to a standard
one), as in~$2x+dx\adequal 2x$.%
\footnote{On occasion Leibniz used the notation ``$\adequal$'' for the
relation of equality.  Note that Leibniz also used our ``$=$'' and
other signs for equality, and did not distinguish between ``$=$'' and
``$\adequal$'' in this regard.  To emphasize the special meaning
\emph{equality} had for Leibniz, it may be helpful to use the symbol
$\adequal$ so as to distinguish Leibniz's equality in a generalized
sense of ``up to'' from the modern notion of equality ``on the
nose.''}

As far as Skolem's nonstandard extension~$\N\subset{}^\ast\N$ is
concerned, anything involving the actual construction of the number
system and the entities called \emph{numbers} would go under the
heading of the \emph{ontology} of mathematical entities.  Note that
\emph{the first order theories} of~$\mathbb{N}$
and~${}^\ast\mathbb{N}$ are identical, as shown by Skolem (for more
details see Section~\ref{FOL}).  In this sense, not only is one not
adding \emph{a lot}, but in fact \emph{one is not adding anything at
all} at the level of the \emph{theory}.

What about the claim that \emph{Leibniz did not have this
perspective}?  It is true that he did not have our perspective on the
\emph{ontological} issues involved in a modern construction of a
suitable number system incorporating infinite numbers, but this
needn't affect the \emph{procedural} match.

What about the claim that \emph{one has to build up a considerable
conceptual machinery to emulate Leibniz's probably rather modest
arsenal of procedural moves; that is to say, we may be able to emulate
all of Leibniz's moves in the modern framework, but it also enables us
to carry out many moves that Leibniz would have never dreamt of}?  As
mentioned above, this is not the case at the level of first order
logic.

What about the claim that \emph{the difference between Leibniz
framework and the infinitesimal framework is large}?  At the
procedural level this is arguably not the case.

What about the claim that \emph{the Weierstrass framework may not
cover all of Leibniz's, but the \emph{difference},
$\text{(A)}-\text{LEI}$, is \emph{small}, indicating that (A) and LEI
are more similar to each other than (B) and LEI are, affecting the
assessment of appropriateness}?  What needs to be pointed out here is
that actually the considerable distance in ontology between (A) and
LEI is about the same is the distance between (B) and LEI.  The
Weierstrassian punctiform continuum where almost all real numbers are
undefinable (so that no individual number of this sort can ever be
specified, unlike~$\pi$,~$e$, etc.) is a far cry from anything one
might have imagined in the 17th century.

As far as the question \emph{Could it be that neither Weierstrass nor
the infinitesimal account is appropriate to cope with Leibniz?}  this
is of course possible in principle.  However, we are interested here
in the practical issue of modern commentators missing some compelling
aspects of interpretation of Leibniz's work because of a self-imposed
limitation to a Weierstrassian interpretive framework.

%Check out what he does for Leibniz on page~30 of his book!  This might
%require some elaboration.

\section{A lid on ontology}

It could be objected that one cannot escape so easily with the general
argument along the lines of ``Let's Ignore (ontological)
Differences,'' or LID for short (putting a \emph{lid} on ontology, so
to speak).

The LID proceeds as follows.  We start with the `real' L, i.e., the
mathematician who lived, wrote, and argued in the 17th century. It
seems plausible to assume that L based his reasoning on a mixture of
first and second order arguments, without clearly differentiating
between the two.

In a reconstruction of~L's arguments, one replaces the cognitive
agent~L by a substitute L1 who argues only in a first order
framework. This entails, in particular, that L1 cannot distinguish
between~$\mathbb{N}$ and~${}^\ast\mathbb{N}$.

However, it seems likely that L could distinguish the two structures,
simply because he did not distinguish between the first and second
order levels.  In other words, the LID recommendation does not help
because the distinction between first and second order does not only
affect the ontology but also the epistemology of the historical agents
involved.

In sum, a modern infinitesimal reconstruction of L deals with a
first-order version of L, namely L1, and not with L.  In line with his
position on geometric algebra, \cite{Un76} could point out that L1 is
a modern artefact, different from the ``real'' L.  Therefore
additional arguments are needed in favor of the hypothesis that L and
L1 are epistemologically sufficiently similar, but this seems
difficult.  In any case, a purely ontological assumption does not
suffice.

To respond to the L \emph{vs} L1 distinction, note that the tools one
needs are \emph{almost} limited to first order logic, but not quite,
since one needs the shadow principle and the relation of infinite
proximity.  Rather than arguing that~$L=L1$, we are arguing that
$L=L1+\epsilon$.

Now the difference between calculus and analysis is that in calculus
one deals mostly with first order phenomena (with the proviso as in
Section~\ref{three}), whereas in analysis one starts tackling
phenomena that are essentially second order, such as the completeness
property i.e., existence of the least upper bound for an arbitrary
bounded set, etc.  It seems reasonable to assume that what they were
doing in the 17th century was calculus rather than analysis.

As far as Unguru is concerned, he is unlikely to be impressed by
interpretations of Leibnizian infinitesimals as quantified
propositions or for that matter by reading Leibniz as if he had
already read not only Weierstrass but also Russell \`a la Ishiguro,
contrary to much textual evidence in Leibniz himself.  We provide a
rebuttal of the Ishiguro--Arthur \emph{logical fiction} reading in
\cite{Ba16} and \cite{Ba16b}.

\section{Robinson, Cassirer, Nelson}

\subsection{Robinson on second-order logic}

The following quote is from Robinson's \emph{Non-standard Analysis}:

\begin{quote}
        The axiomatic systems for many algebraic concepts such as
groups or fields are formulated in a natural way within a first order
language\ldots, However, interesting parts of the \emph{theory} of
such a concept may well extend beyond the resources of a first order
language. Thus, in the theory of groups statements regarding
subgroups, or regarding the existence of subgroups of certain types
will, in general, involve quantification with respect to sets of
individuals\ldots{} \cite[section~2.6, p.~19]{Ro66}
\end{quote}

This appears to amount to a claim that the local ontology may indeed
be often formulated in first-order terms, while the global ontology is
deeply infected by second-order concepts.  The latter may typically
involve objects and arguments qualified as second or higher order with
respect to the former, which nevertheless are of the first-order type
when considered as related to the background set universe.

\vfill\eject

\subsection{Ernst Cassirer}

Does the equation~$L = L1 +\epsilon$ not amount to an underestimation
of the historical Leibniz?  Is it reasonable to assume that he only
invented the calculus, and not analysis?  According to Cassirer, the
basic concepts of analysis were deeply soaked with philosophy, i.e.,
for Leibniz mathematical and philosophical concepts were intimately
related:%
\footnote{In support of this claim, Cassirer refers here in particular
to Gottfried Wilhelm Leibniz, Die philosophischen Schriften, hrg. von
Carl Immanuel Gerhardt, 7 Bde., Berlin 1875--1890, Bd. VII,
S. 542. (Cassirer 1902, p. xi)}
\begin{quote}
Leibniz himself asserted that the new analysis has sprung from the
innermost source of philosophy, and he assigned to both regions [i.e.,
analysis and philosophy] the task to confirm and to elucidate each
other.%
\footnote{In the original: ``Leibniz selbst hat es ausgesprochen,
da\ss{} die neue Analysis aus dem innersten Quell der Philosophie
geflossen ist, und beiden Gebieten die Aufgabe zugewiesen, sich
wechselseitig zu best\"atigen und zu erhellen.''}
\cite[p.~xi]{Ca02}
\end{quote}

If L = L1 + epsilon, i.e., if the historical Leibniz was mainly
dealing with \emph{calculus}, this may appear hardly compatible with
Cassirer's perspective; see \cite{MK}.  This impression would be,
however, a misunderstanding. In order to forestall it, it merits being
pointed out that developing the calculus was a great mathematical
achievement of philosophical relevance.  It is only today that the
term \emph{calculus} possesses a connotation of routine undergraduate
mathematics, but not in the 19th century.

As far as Cassirer is concerned, Leibniz was indeed doing
\emph{analysis} as l'H\^opital called it.  It is not even sure
Cassirer was aware of the more advanced analysis.  Leibnizian calculus
only seems ``trivial" from the standpoint of properly 20th century
mathematics.  It is an advance in understanding when we make a
distinction between Leibnizian calculus and analysis.  We don't mean
to diminish Leibniz's greatness by this distinction, nor do we suggest
that Cassirer was \emph{wrong}.  He was merely using the term
\emph{analysis} in its 17--19th century sense rather than the sense in
which we use it today.

\subsection{Edward Nelson}

As far as the passage from \cite{Ro66} is concerned, we find the
following comment at the end of the paragraph:
\begin{quote}
The following framework for higher order structures and higher order
languages copes with these and similar cases. It is rather
straightforward and suitable for our purposes. \cite[section~2.6,
p.~19]{Ro66}
\end{quote}
Robinson then proceeds to develop a solution, which roughly
corresponds to level~(3) as outlined in Section~\ref{FOL}.  The claim
that we are dealing with first order logic plus \emph{standard part}
is in a sense a mathematical theorem, undermining the contention that
``this seems hardly compatible, etc.''

Edward Nelson demonstrated that infinitesimals can be found within the
ordinary real line itself in the following sense. Nelson finds
infinitesimals in the real line by means of enriching the language
through the introduction of a unary predicate \emph{standard} and an
axiom schemata (of Idealization), one of most immediate instances of
which implies the existence of infinitely large integers and hence
nonzero infinitesimals; see \cite{Ne77}.  This is closely parallel to
the dichotomy of \emph{assignable} vs \emph{inassignable} in Leibniz,
whereas Carnot spoke of \emph{quantit\'es d\'esign\'ees} \cite{Ca97},
\cite[p.~46]{Ba89}.  Thus, we obtain infinitesimals as soon as we
assume that (1) there are assignable (or standard) reals, that obey
the same rules as all the reals, and (2) there are reals that not
assignable.

In more detail, one considers the ordinary ZFC formulated in first
order logic (here the term is used in a different sense from the rest
of this article), adds to it the unary predicate and the axiom
schemata, and obtains a framework where calculus and analysis can be
done with infinitesimals.  For further discussion see \cite{KK15}.

The passage from Robinson cited above does indicate that second order
theory may often be interesting.  However, in the case of the
calculus/analysis as it was practiced in the 17th century, we are not
aware of a single significant result that cannot be formulated in a
system of type $\text{FOL}+\square$ (see Section~\ref{one} for
examples of results that can).  Arguably it was calculus (rather than
analysis) that Leibniz invented, in the sense that there don't appear
to be any essentially second order statements there.

It may seem surprising that there could be a kind of pre-established
harmony between a modern logical category, namely, first-order
results, and a historical category, namely, results of 17th century
calculus.  This idea suggests further questions: does this only hold
for the calculus, or is it always or often the case that a
historically earlier realization of a theory covers only the
first-order part of its successor.  How do arithmetic and geometry
behave in this respect?  Would it really make sense to systematically
distinguish between Leibniz and Leibniz1, Euclid and Euclid1, etc.?

\subsection{Not standalone}

Let us return to the comment \emph{In the context of Skolem's
construction of nonstandard natural numbers and some related stuff,
one is impressed by how many nonstandard numbers one needs to move
from~$\mathbb N$ to~${}^\ast\mathbb{N}$, e.g., in Henkin's model
${}^\ast\mathbb N = \mathbb N + (\Z \times \Q)$ and it goes without
saying that Leibniz did not have this perspective}, that was addressed
briefly above.  One could elaborate on the ``impression'' concerning
``how many nonstandard numbers'' one needs to
define~${}^\ast\mathbb{N}$ consistently and conveniently.

An infinitely large number say~$H$ is not a standalone object, but
rather lives in a community of numbers obeying certain laws which
mathematicians anticipate as a goal of the construction of a
nonstandard number system~${}^\ast\mathbb N$.  Such a commitment to
anticipated laws forces Skolem and others to add to~$\mathbb N$ a
suitable entourage of~$H$ along with~$H$ itself.  What are the laws
involved?

Modern specialists in Nonstandard Analysis (NSA) stipulate
that~${}^\ast \mathbb N$ should satisfy the axioms of Peano Arithmetic
and moreover, satisfy the same sentences of the language (not
necessarily consequences of the axioms) that are true in~$\mathbb N$
itself.  This is called (the principle of) Transfer today.
Mathematicians of the 17th (or even 19th) century had neither this
perspective nor the tools consistently to define~${}^\ast\mathbb N$
or~${}^\ast\mathbb R$.

On the other hand, one can argue that there is no need for actually
rigorously defining~${}^\ast\mathbb N$ in order to make use of its
benefits.  One can argue that it is sufficient to have some idea of
Transfer on top of an acceptance of infinitely large numbers per se
(possibly as \emph{useful fictions}, to borrow Leibniz's expression).
We have argued that the Leibnizian \emph{Law of continuity} is closely
related to the Transfer principle; see \cite{KS1}.

Therefore the claim that \emph{Leibniz had not the slightest idea of
this stuff} (the ``stuff'' being the modern technique of building
nonstandard models) is perhaps technically true, but it does not
reflect all the aspects of the interrelations within the
Leibniz/Weierstrass/NSA triangle.

\section{Conclusion}

The vast oeuvre of Leibniz is still in the process of publication.  In
principle a lucky scholar might one day unearth a manuscript where
Leibniz tackles a property equivalent to the completeness of the reals
(after all the existence of the shadow is so equivalent), involving
quantification over \emph{all} sets of the number system and therefore
second-order.

However, this is unlikely in view of Leibniz's reluctance to deal with
infinite collections, as mentioned above.  If level (3) of first order
logic is indeed suitable for developing modern proxies for the
inferential moves found in Leibnizian infinitesimal calculus, as we
have argued, then modern infinitesimal frameworks are more appropriate
to interpreting Leibnizian infinitesimal calculus than modern
Weierstrassian ones.

\section*{Acknowledgments}

M.~Katz was partially funded by the Israel Science Foundation grant
number~1517/12.

\bigskip\noindent \textbf{Piotr B\l aszczyk} is Professor at the
Institute of Mathematics, Pedagogical University (Cracow, Poland). He
obtained degrees in mathematics (1986) and philosophy (1994) from
Jagiellonian University (Cracow, Poland), and a Ph.D. in ontology
(2002) from Jagiellonian University.  He authored \emph{Philosophical
Analysis of Richard Dedekind's memoir Stetigkeit und irrationale
Zahlen} (Cracow, 2008, Habilitationsschrift).  He co-authored
\emph{Euclid, Elements, Books V--VI.  Translation and commentary}
(Cracow, 2013), and \emph{Descartes, Geometry. Translation and
commentary} (Cracow, 2015).  His research interest is in the idea of
continuum and continuity from Euclid to modern times.

\medskip\noindent \textbf{Vladimir Kanovei} graduated in 1973 from
Moscow State University, and obtained a Ph.D. in physics and
mathematics from Moscow State University in 1976. In 1986, he became
Doctor of Science in physics and mathematics at Moscow Steklov
Mathematical Institute (MIAN).  He is currently Principal Researcher
at the Institute for Information Transmission Problems (IPPI), Moscow,
Russia, and Professor at the Moscow State University of Railway
Engineering (MIIT), Moscow, Russia.  Among his publications is the
book \emph{Borel equivalence relations. Structure and classification},
University Lecture Series 44, American Mathematical Society,
Providence, RI, 2008.

\bigskip\noindent \textbf{Karin Katz} (B.A. Bryn Mawr College, '82);
Ph.D. Indiana University, '91) teaches mathematics at Bar Ilan
University, Ramat Gan, Israel.  Among her publications is the joint
article ``Proofs and retributions, or: why Sarah can't \emph{take}
limits'' published in \emph{Foundations of Science}.

\medskip\noindent \textbf{Mikhail G. Katz} (B.A. Harvard University,
'80; Ph.D. Columbia University, '84) is Professor of Mathematics at
Bar Ilan University.  Among his publications is the book
\emph{Systolic geometry and topology}, with an appendix by Jake
P. Solomon, published by the American Mathematical Society; and the
article (with T.~Nowik) ``Differential geometry via infinitesimal
displacements'' published in \emph{Journal of Logic and Analysis}.

\medskip\noindent \textbf{Taras S. Kudryk} (born 1961, Lviv, Ukraine)
is a Ukrainian mathematician and associate professor of mathematics at
Lviv National University.  He graduated in 1983 from Lviv University
and obtained a Ph.D. in physics and mathematics in 1989.  His main
interests are nonstandard analysis and its applications to functional
analysis. He is the author of books about nonstandard analysis (in
Ukrainian and English) and textbooks about functional analysis (in
Ukrainian) co-authored with V. Lyantse.  Kudryk has performed research
in nonstandard analysis in collaboration with V. Lyantse and Vitor
Neves.  His publications appeared in \emph{Matematychni Studii},
\emph{Siberian Journal of Mathematics}, and \emph{Logica Universalis}.

\medskip\noindent \textbf{Thomas Mormann} is Professor at the
Departamento de L\'ogica y Filosof\'\i a de la Ciencia de la
Universidad del Pa\'\i s Vasco UPV/EHU (Donostia-San Sebasti\'an,
Spain).  He obtained a PhD in Mathematics from the University of
Dortmund, and habilitated in Philosophy, Logic, and Philosophy of
Science at the University of Munich.  He published numerous papers in
philosophy of science, history of philosophy of science, epistemology,
and related areas.  He is the editor of Carnap's ``Anti-Metaphysical
Writings" (in German).  Presently, his main research interest is in
the philosophy of Ernst Cassirer and, more generally, in the Marburg
Neo-Kantianism.

\medskip\noindent \textbf{David Sherry} is Professor of Philosophy at
Northern Arizona University, Flagstaff, AZ.  E-mail:
david.sherry@nau.edu.  He has research interests in philosophy of
mathematics, especially applied mathematics and nonstandard analysis.

\end{document}